\documentclass[12pt]{amsart}
\usepackage{latexsym,amssymb,amscd}
\def\frk{\frak}               

\def\mm{{\frk m}}

\def\B'c{{\mathcal{B'}}}
\def\U'c{{\mathcal{U'}}}

\def\eb{{\bold e}}
\def\lb{{\bold l}}

\def\xb{{\bold x}}

\def\opn#1#2{\def#1{\operatorname{#2}}} 
\opn\chara{char}
\opn\length{\ell}
\opn\projdim{proj\,dim}
\opn\injdim{inj\,dim}
\opn\ini{{\rm in}}
\opn\rank{rank}
\opn\Tiefe{Tiefe}
\opn\grade{grade}
\opn\depth{depth}
\opn\height{height}
\opn\embdim{emb\,dim}
\opn\codim{codim}

\opn\Tr{Tr}
\opn\bigrank{big\,rank}
\opn\superheight{superheight}\opn\lcm{lcm}
\opn\trdeg{tr\,deg}%
\opn\reg{\rm reg}
\opn\lreg{lreg}
\opn\hdeg{{\rm hdeg}}
\opn\ideg{ideg}
\opn\pd{{\rm pd}}
\opn\div{div}
\opn\Div{Div}
\opn\cl{cl}
\opn\Cl{Cl}
\opn\Spec{Spec}
\opn\Supp{Supp}
\opn\supp{supp}
\opn\Sing{Sing}
\opn\Ass{Ass}
\opn\Ann{Ann}
\opn\Rad{Rad}
\opn\Soc{Soc}
\opn\Ker{Ker}
\opn\Coker{Coker}
\opn\Im{Im}
\opn\Hom{Hom}
\opn\Tor{Tor}
\opn\Ext{Ext}
\opn\End{End}
\opn\Aut{Aut}
\opn\id{id}

\opn\nat{nat}
\opn\GL{GL}
\opn\SL{SL}
\opn\mod{mod}
\opn\ord{ord}
\opn\Gin{Gin}
\opn\aff{aff}
\opn\con{conv}
\opn\relint{relint}
\opn\st{st}
\opn\lk{lk}
\opn\cn{cn}
\opn\core{core}
\opn\vol{vol}
\opn\gr{gr}

\def\pot#1#2{#1[\kern-0.28ex[#2]\kern-0.28ex]}

\opn\dirlim{\underrightarrow{\lim}}
\opn\invlim{\underleftarrow{\lim}}
\let\union=\cup

\let\dirsum=\oplus
\let\tensor=\otimes
\let\iso=\cong

\let\Dirsum=\bigoplus

\def\pnt{{\raise0.5mm\hbox{\large\bf.}}}
\def\lpnt{{\hbox{\large\bf.}}}
\let\to=\rightarrow

\def\Implies{\ifmmode\Longrightarrow \else
     \unskip${}\Longrightarrow{}$\ignorespaces\fi}
\def\implies{\ifmmode\Rightarrow \else
     \unskip${}\Rightarrow{}$\ignorespaces\fi}
\def\iff{\ifmmode\Longleftrightarrow \else
     \unskip${}\Longleftrightarrow{}$\ignorespaces\fi}

\let\:=\colon
\newtheorem{Theorem}{Theorem}[section]
\newtheorem{Lemma}[Theorem]{Lemma}
\newtheorem{Corollary}[Theorem]{Corollary}
\newtheorem{Proposition}[Theorem]{Proposition}

\newtheorem{Construction}[Theorem]{Construction}
\let\epsilon=\varepsilon
\let\phi=\varphi
\let\kappa=\varkappa
\title{Bounds for Betti numbers}
\author{Tim R\"omer}

\begin{document}

\maketitle
\begin{abstract}
In this paper we prove parts of a conjecture of Herzog
giving lower bounds on the rank
of the free modules appearing in the linear strand
of a graded $k$-th syzygy module 
over the polynomial ring. If in addition the module
is $\mathbb{Z}^n$-graded we show that
the conjecture holds in full generality.
Furthermore, 
we give lower and upper bounds for the 
graded Betti numbers
of graded ideals with a linear resolution and a fixed number of 
generators.
\end{abstract}

\section*{Introduction}
Let $S=K[x_1,\ldots,x_n]$ be the polynomial ring over a field $K$
equipped with the standard grading by setting $\deg(x_i)=1$, and 
let $M$ be a finitely generated graded $S$-module. 
We denote by $\beta_{i,i+j}(M)=\dim_K \Tor_i(M,K)_{i+j}$ the graded Betti numbers of $M$.

Assume that the initial degree of $M$ is $d$, i.e.\ we have $M_i=0$ for $i<d$ 
and $M_d \neq 0$. 
We are interested in the 
numbers $\beta_i^{lin}(M)=\beta_{i,i+d}(M)$ for 
$i \geq 0$. These numbers determine 
the rank of the free modules appearing in the linear strand of the  
minimal graded free resolution of $M$. 
Let $p=\max\{i \colon \beta_i^{lin}(M)\neq 0 \}$ be the length of the 
linear strand. In \cite{HE} Herzog conjectured the following:
\vspace{3pt}
\\
{
\label{herzog} 
\bf
Conjecture. 
\it
Let $M$ be a $k$-th syzygy module whose linear strand has length $p$, then
$$
\beta^{lin}_{i}(M)\geq \binom{p+k}{i+k}
$$ 
for $i=0,\ldots,p$.
}
\vspace{3pt}

This conjecture is motivated by a result of Green \cite{GR} 
(see also Eisenbud and Koh \cite{EIKO})
that contains the case $i=0, k=1$. 
For $k=0$ these lower bounds were shown by Herzog \cite{HE},
and Reiner, Welker \cite{REWE} proved them for $k=1$ if $M$ is a monomial ideal.

In this paper we prove the conjecture for $k=1$.  
For $k>1$ we get the weaker result:  
$$
\text{If }\beta^{lin}_{p}(M) \neq 0  \text{ for $p>0$ and $M$ is a 
$k$-th syzygy module, then } \beta^{lin}_{p-1}(M) \geq p+k. 
$$
We also show that the conjecture holds in full generality for
finitely generated $\mathbb{Z}^n$-graded $S$-modules. The first three 
sections of this paper are concerned with the question above. 

In the last years many authors (see for example \cite{BI,HU,PA}) were interested 
in the following problem: 
Fix a possible Hilbert function $H$ for a graded ideal.
Let $\mathcal{B}(H)$ be the set of Betti sequences $\{\beta_{i,j}(I)\}$ where $I \subset S$ is a
graded ideal with Hilbert function $H$.
On $\mathcal{B}(H)$ we consider a partial 
order: We set $\{\beta_{i,j}(I)\}\geq \{\beta_{i,j}(J)\}$ 
if $\beta_{i,j}(I)\geq \beta_{i,j}(J)$ for all $i,j \in \mathbb{N}$.
It is known that $\mathcal{B}(H)$ has a unique maximal element given
by the Betti sequence of the lex-segment ideal in the family of considered ideals.
In general there is more than one minimal element (see \cite{CHEV}). 

In Section $4$ we study a related problem. We fix an integer $d \geq 0$ and
$0 \leq k \leq\binom{n+d-1}{d}$.
Let $\mathcal{B}(d,k)$ be the set of Betti sequences $\{\beta_{i,j}(I)\}$ where 
$I \subset S$ is a
graded ideal with $d$-linear resolution and $\beta_{0,d}(I)=k$. 
We show that, independent of the characteristic of the base field,  
there is a unique minimal and a unique maximal element in $\mathcal{B}(d,k)$.

The author is grateful to Prof.\ Herzog for inspiring discussions on the 
subject of the paper.

\section{Preliminaries on Koszul complexes}
Let $K$ be a field, $V$ be an $n$-dimensional $K$-vector space 
with basis $\xb=x_1,\ldots,x_n$ and
$S=K[V]$ be the symmetric algebra over $V$ equipped with 
the standard grading by setting $\deg(x_i)=1$. 
Furthermore let $\mm=(x_1,\ldots,x_n)$ be the graded maximal ideal of $S$ and
$0 \neq M$ be a finitely generated graded $S$-module which is generated 
in not negative degrees, i.e.\ $M_i=0$ for $i<0$. 

Consider a graded free $S$-module $L$ of rank $j$ which 
is generated in degree $1$
and let $\bigwedge L$ be the exterior algebra over $L$.
Then $\bigwedge L$ inherits the structure of a bigraded $S$-module. 
If $z \in \bigwedge^i L$ and $z$ has $S$-degree $k$,
then we give $z$ the bidegree $(i,k)$.
We call $i$ the homological degree (hdeg for short) and $k$ the
internal degree (deg for short) of $z$. 
 
We consider maps $\mu \in L^*=\Hom_S(L,S)$. 
Note that $L^*$ is again a graded free $S$-module generated in degree $-1$.
It is well-known (see \cite{BRHE}) that  
$\mu$ defines a graded $S$-homomorphism $\partial_{\mu}\colon \bigwedge L \to \bigwedge L$ 
of (homological) degree $-1$.

Recall that if we fix a basis $e_1,\ldots,e_j$ of $L$,
then $\bigwedge^i L$ is the graded free $S$-module with 
basis 
consisting of all monomials
$e_J=e_{j_1}\wedge\ldots\wedge e_{j_i}$ 
with $J=\{j_1<\ldots<j_i \}\subseteq [j]=\{1,\ldots,j \}$. 
One has
$$\partial_{\mu}(e_{j_1}\wedge\ldots\wedge e_{j_i})=
\sum_{k=1}^{i}(-1)^{\alpha(k,J)}\mu(e_{j_k}) e_{j_1}\wedge\ldots 
\hat{e}_{j_k}\ldots\wedge e_{j_i}$$ 
where 
for $F,G \subseteq [n]$ we set $\alpha(F,G)=|\{(f,g) \colon f>g,f \in F,g \in G \}|$
and where $\hat{e}_{j_k}$ indicates that $e_{j_k}$ is to be omitted from 
the exterior product.
Denote by $e_1^*,\ldots,e_j^*$ the basis of $L^*$ with $e_i^*(e_i)=1$ 
and $e_i^*(e_k)=0$ for $k \neq i$.
In order to simplify notations we set $\partial_{i}=\partial_{e_i^*}$. 
Then $\partial_{\mu}=\sum_{k=1}^j \partial_{\mu}(e_k) \partial_k$, and 
we have the following:
\begin{Lemma}
\label{rules}
\it
Let $z,\Tilde{z} \in \bigwedge  L$ be bihomogeneous elements,
$f \in S$ and $\mu,\nu \in L^*$. Then:
\begin{enumerate}
\renewcommand{\labelenumi}{(\roman{enumi})}
\item 
$f \partial_{\mu}=\partial_{f \mu}$,
\item 
$\partial_{\mu}+\partial_{\nu}=\partial_{\mu+\nu}$,
\item 
$\partial_{\mu}\circ \partial_{\mu}=0$,
\item 
$\partial_{\mu}\circ\partial_{\nu}=-\partial_{\nu}\circ\partial_{\mu}$,
\item 
$\partial_{\mu}(z \wedge \Tilde{z})=\partial_{\mu}(z)\wedge 
\Tilde{z}+(-1)^{\hdeg(z)}z \wedge \partial_{\mu}(\Tilde{z})$.
\end{enumerate}
\end{Lemma}
\begin{proof}
Straightforward calculations (most of them are done in \cite{BRHE}).
\end{proof} 

We fix a graded free $S$-module $L$ of rank $n$ for the rest of the paper.
Let $\eb=e_1,\ldots,e_n$ be a basis of $L$ with $\deg(e_i)=1$ for $i \in [n]$, 
and $\mu \in L^*$ with 
$\mu(e_i)=x_i$ for $i \in [n]$.
For $j=1,\ldots,n$ let $L(j)$ be the graded free submodule of $L$
generated by $e_1,\ldots,e_j$.
Then $(K(j),\partial)$ is the Koszul complex of $x_1,\ldots,x_j$ with 
values in $M$ where $K(j)=\bigwedge L(j)\tensor_S M$ and $\partial$ is the 
restriction of $\partial_{\mu}\tensor_S \id_M$ to $\bigwedge L(j)\tensor_S M$.

We denote by $H(j)$ the homology of the complex $K(j)$, and the 
homology class of a cycle $z \in K(j)$ will be denoted by $[z]$.

Notice that $K_i(j)_{i+k}=0$ for $k<0$ and that 
$H_i(n)\iso\Tor_i(K,M)$ are isomorphic as graded $K$-vector spaces.
One has the following exact sequence (see \cite{BRHE}):
$$
\ldots \to H_{i+1}(j)\to H_{i}(j-1)(-1) \to H_{i}(j-1)\to H_{i}(j)\to \ldots.
$$

The following observation is crucial for the rest of the paper.
For a homogeneous element $z \in K_i(j)$ we can 
write $z$ uniquely as
$z=e_k \wedge\partial_k(z) + r_z$ and $e_k$ divides none of the monomials
of $r_z$.
\begin{Lemma}
\label{cycle}
\it
Let $z \in K_{i}(j)$ be a homogeneous cycle of bidegree $(i,l)$.
Then $\partial_k(z)$ is for all $k \in [n]$ a homogeneous cycle of bidegree $(i-1,l-1)$. 
\end{Lemma}
\begin{proof}
This follows from \ref{rules}.
\end{proof}

In the sequel we need the following:
\begin{Lemma}
\label{homology}
\it
Let $p \in \{0,\ldots,n\}$ and $t \in \mathbb{N}$.
Suppose that
$H_p(j)_{p+l}=0$ for $l=-1,\ldots, t-1$. 
Then:
\begin{enumerate}
\renewcommand{\labelenumi}{(\roman{enumi})}
\item
$H_p(j-1)_{p+l}=0$ for $l=-1,\ldots, t-1$, 
\item
$H_p(j-1)_{p+t}$ is isomorphic to a submodule of $H_p(j)_{p+t}$, 
\item
$H_i(j)_{i+l}=0$ for $l=-1,\ldots, t-1$ and $i=p, \ldots, j$. 
\end{enumerate}
\end{Lemma}
\begin{proof}
We prove (i) by induction on $l \in \{-1, \ldots, t-1 \}$. If $l=-1$
there is nothing to show because $H_p(j-1)_{p+l}=0$ for $l<0$. 
Now let $l>-1$ and consider the exact sequence 
$$
\ldots \to H_{p}(j-1)_{p+l-1} \to H_{p}(j-1)_{p+l}\to H_{p}(j)_{p+l}\to \ldots.
$$
Since by induction hypothesis $H_{p}(j-1)_{p+l-1}=0$ and by the assumption 
$H_{p}(j)_{p+l}=0$, we get that $H_{p}(j-1)_{p+l}=0$.

For $l=t$ the exact sequence of the Koszul homology together with (i) yields
$$
0 \to H_{p}(j-1)_{p+t}\to H_{p}(j)_{p+t}\to \ldots,
$$
which proves (ii).

We show (iii) by induction on $j \in [n]$. The case $j=1$ is trivial and
for $j>1$ and $i=p$ the assertion is true by assumption.
Now let $j>1$, $i>p$ and consider
$$
\ldots \to H_{i}(j-1)_{i+l} \to H_{i}(j)_{i+l}\to H_{i-1}(j-1)_{i-1+l}\to \ldots.
$$
By (i) and the induction on $j$ we get that
$H_{i}(j-1)_{i+l}=H_{i-1}(j-1)_{i-1+l}=0$. 
Hence $H_{i}(j)_{i+l}=0$.
\end{proof}

\section{Lower bounds for Betti numbers of graded $S$-modules}
In this section 
$M$ is always a finitely generated graded $S$-module which is generated 
in degrees $\geq 0$. 
For $0 \neq z \in K_i(j)$
we write 
$$
z=m_J e_J+\sum_{I \subseteq [n], I \neq J} m_I e_I
$$
with coefficients in
$M$, and where $e_J$ is the lexicographic largest 
monomial of all $e_L$ with $m_L \neq 0$. Recall
that for $I,J \subseteq [n]$, 
$I=\{i_1<\ldots< i_t \}$, $J=\{j_1<\ldots< j_{t'} \}$ 
$e_I <_{lex} e_J$ if either $t<t'$ or $t=t'$ and there exists a number $p$ with 
$i_l=j_l$ for $l<p$ and $i_p>j_p$.
We call $\ini (z)=m_J e_J$ the initial term of $z$.
Furthermore for $I=\{i_1<\ldots< i_t \}\subseteq [n]$ we write
$\partial_I=\partial_{i_1}\circ\ldots\circ\partial_{i_t}$.

\begin{Lemma}
\label{linear}
Let $p \in \{0,\ldots,j\}$, $r \in \{0,\ldots,p \}$ and
$0 \neq z \in K_p(j)$ be homogeneous with $\ini(z)=m_J e_J$. 
Then for all $I \subseteq J$ with $|I|=r$ the elements $\partial_I(z)$ are 
$K$-linearly independent in $K_{p-r}(j)$.
In particular, if $z$ is a cycle, then 
$\{\partial_I(z) \colon I \subseteq J, |I|=r \}$ is a set of $K$-linearly independent 
cycles.
\end{Lemma}
\begin{proof}
This follows from the fact that $\ini(\partial_I(z))=m_J e_{J-I}$.
Induction on $r \in \{0,\ldots,p \}$ proves that all $\partial_I(z)$ are cycles if $z$ is one.
\end{proof}

\begin{Lemma}
\label{wichtig}
\it
Let $p \in [j]$, $t \in \mathbb{N}$ and
$z \in K_p(j)_{p+t}$. 
Assume that $H_{p-1}(j)_{p-1+l}=0$ for $l=-1,\ldots,t-1$.
\begin{enumerate}
\renewcommand{\labelenumi}{(\roman{enumi})}
\item
If $p<j$ and $\partial_j(z)=\partial(y)$ for some $y$, then there exists $\Tilde{z}$ such that
$\Tilde{z}=z+\partial(r)$ and $\partial_j(\Tilde{z})=0$.
In particular, $\Tilde{z} \in K_p(j-1)$, and if $z$ is a cycle, then $[z]=[\Tilde{z}]$.
\item
If $p=j$ and $\partial_j(z)=\partial(y)$ for some $y$, then $z=0$.
In particular, if $z \neq 0$ is a cycle, then we always 
have $0 \neq [\partial_j(z)] \in H_{p-1}(j)_{p-1+t}$.
\end{enumerate}
\end{Lemma}
\begin{proof}
We proceed by induction on $t \in \mathbb{N}$ to prove (i). If $t=0$, then 
$y \in K_{p}(j)_{p+t-1}=0$, and so 
$\partial_j(z)=0$. Thus we choose $\Tilde{z}=z$. 

Let $t>0$
and assume that $\partial_j(z)=\partial(y)$. We see
that $\partial_j(y)$ is a cycle because
$$
0=\partial_j (\partial_j(z))=\partial_j (\partial(y))=- \partial(\partial_j(y)).
$$
But $\partial_j(y) \in K_{p-1}(j)_{p-1+t-1}$. Since $H_{p-1}(j)_{p-1+t-1}=0$, it
follows that $\partial_j(y)=\partial(y')$ is a boundary for some element $y'$.

By the induction hypothesis we get $\Tilde{y}=y+\partial(r')$
such that $\partial_j(\Tilde{y})=0$. Note  that
$$
\partial(\Tilde{y})=\partial(y)=\partial_j(z).
$$
We define 
$$ 
\Tilde{z}=z + \partial(e_j \wedge\Tilde{y} )=z +x_j \Tilde{y}- e_j \wedge\partial_j(z).
$$
Then
$$
\partial_j(\Tilde{z})=\partial_j(z)+x_j \partial_j(\Tilde{y})-\partial_j(e_j)\wedge\partial_j(z) 
+e_j \wedge \partial_j \circ \partial_j(z)
=\partial_j(z)-\partial_j(z)=0
$$
and this proves (i).

If $p=j$, we see that $z=m e_{[j]}$ for some $m \in M$ and therefore
$\partial_j(z)\neq 0$ if and only if $z \neq 0$.

We prove (ii) by induction on $t \in \mathbb{N}$. For $t=0$ there is nothing to show.
Let $t>0$ and assume $\partial_j(z)=\partial(y)$.
By the same argument as in the proof of (i) we get $\partial_j(y)=\partial(y')$ 
for some $y'$.
The induction hypothesis implies $y=0$, and then $z=0$.
\end{proof}

\begin{Lemma}
\label{wichtig2}
\it
Let $p \in \{0,\ldots,j\}$, $t \in \mathbb{N}$
and
$0 \neq z \in K_p(j)_{p+t}$. 
Assume that $H_{p}(j)_{p+l}=0$ for $l=-1,\ldots,t-1$ and let $q \in [j]$.
If $z \in K_p(j-q)_{p+t}\subseteq K_p(j)_{p+t}$ and $\partial(y)=z$ in $K(j)$ for some element $y$,
then there exists $\Tilde{y}=y+\partial(r) \in K_{p+1}(j-q)_{p+1+t-1}$ such that
$\partial(\Tilde{y})=z$.
\end{Lemma}
\begin{proof}
We prove the assertion by induction on $j-q$. For $j=q$ there is nothing to show. Let $j>q$.
Since $0=\partial_j(z)$
and
$$
z=\partial(y)=\partial(e_j \wedge\partial_j(y)+r)=
x_j \partial_j(y)
- e_j \wedge \partial(\partial_j(y))
+ \partial(r),
$$
we see that $\partial_j(y)$ is a cycle and therefore a boundary by the assumption
that $H_p(j)_{p+t-1}=0$.
By \ref{wichtig} we 
may assume that $y \in K(j-1)$.
By the induction hypothesis we find
the desired $\Tilde{y}$ in $K(j-q)$.
\end{proof}

\begin{Lemma}
\label{erstens}
\it
Let $t \in \mathbb{N}$.
If $\beta_{n-1,n-1+l}(M)=0$ for $l=-1,\ldots,t-1$ and $\beta_{n,n+t}(M)\neq 0$,
then there exists a basis $\eb$ of $L$ and
a cycle $z \in K_n(n)_{n+t}$ such that
\begin{enumerate}  
\renewcommand{\labelenumi}{(\roman{enumi})}
\item
$[z]\in H_n(n)_{n+t}$ is not zero,
\item
$[\partial_i(z)] \in H_{n-1}(n)_{n-1+t}$  are $K$-linearly independent for $i=1,\ldots,n$. 
\end{enumerate}
In particular, $\beta_{n-1,n-1+t}(M)\geq n$.
\end{Lemma}
\begin{proof}
Let $\eb$ be an arbitrary basis of $L$. Since $\beta_{n,n+t}(M)\neq 0$ there 
exists a cycle $z \in K_n(n)_{n+t}$ with $0 \neq [z]\in H_n(n)_{n+t}$.
Furthermore $H_{n-1}(n)_{n-1+l}=0$ for $l=-1,\ldots,t-1$.
In this situation we have $z=m e_{[n]}$ for some socle element $m$ of $M$ and
we want to show that every equation 
$$
0 = \sum_{i=1}^n \mu_i [\partial_i(z)]=[\sum_{i=1}^n \mu_i \partial_i(z)]
\quad \text{with} \quad \mu_i \in K,
$$
implies $\mu_i=0$ for all $i \in [n]$.
Assume there is such an equation where not all $\mu_i$ are zero. 
After a base change we may assume that $\sum_{i=1}^n \mu_i \partial_i=\partial_n$.
We get 
$$
0= [\partial_n (z)],
$$ 
contradicting to \ref{wichtig} (ii).
\end{proof}

\begin{Theorem}
\label{Schritt1}
\it
Let $t \in \mathbb{N}$ and $p \in [n]$.
If $\beta_{p-1,p-1+l}(M)=0$ for $l=-1,\ldots,t-1$ and $\beta_{p,p+t}(M)\neq 0$,
then there exists a basis $\eb$ of $L$ and
a cycle $z \in K_p(n)_{p+t}$ such that
\begin{enumerate}  
\renewcommand{\labelenumi}{(\roman{enumi})}
\item
$[z]\in H_p(n)_{p+t}$ is not zero,
\item
$[\partial_i(z)] \in H_{p-1}(n)_{p-1+t}$ are $K$-linearly independent for $i=1,\ldots,p$. 
\end{enumerate}
In particular, $\beta_{p-1,p-1+t}(M)\geq p$.
\end{Theorem}
\begin{proof}
We have $H_p(n)_{p+t}\neq 0$ because $\beta_{p,p+t}(M)\neq 0$. Choose 
$0 \neq h \in H_p(n)_{p+t}$.
We prove by induction on $n$ that we can find a basis $\eb$ of $L(n)$
and a cycle $z \in K_p(n)_{p+t}$ representing $h$ such that every equation
$$
0 = \sum_{i=1}^p \mu_i [\partial_i(z)]=[\sum_{i=1}^p \mu_i \partial_i(z)]
\quad \text{with} \quad \mu_i \in K,
$$
implies $\mu_i=0$ for all $i$. 
The cases $n=1$ and 
$n>1,$ $p=n$ were shown in \ref{erstens}. 

Let $n>1$ and $p<n$.
Assume that there is a basis $\eb$ and such an equation for a cycle $z$ with
$[z]=h$ where not all $\mu_i$ are zero.
After a base change of $L(n)$ we may assume that $\sum_{i=1}^p \mu_i \partial_i=\partial_n$.
Then
$0= [\partial_n (z)],$ 
and therefore $\partial_n (z)=\partial(y)$ for some element $y$. 
By \ref{wichtig} 
we can find an element   
$\Tilde{y}$ such that $[\Tilde{y}]=[z]$ and $\Tilde{y} \in K_p(n-1)_{p+t}$.
Now \ref{homology} guarantees that we can apply our
induction hypothesis to $\Tilde{y}$ and we find a base change $\lb=l_1,\ldots, l_{n-1}$
of $e_1,\ldots,e_{n-1}$, $[\Tilde{z}]=[\Tilde{y}]$ in $H_p(n-1)_{p+t}$ 
(with respect to the new basis) such that
$[\partial_i(\Tilde{z})] \in H_{p-1}(n-1)_{p-1+t}$ are 
$K$-linearly independent for $i=1,\ldots,p$. 
By \ref{homology} we have $H_i(n-1)_{i+t}\subseteq H_i(n)_{i+t}$ for $i=p-1, p$. 
Then $\Tilde{z}$ is the desired cycle because $[\Tilde{z}]=[z]$ in $H_p(n)_{p+t}$.
\end{proof}

The Castelnuovo-Mumford regularity for a finitely generated graded $S$-module $0 \neq M$ 
is defined as
$\reg(M)=\max\{j \in \mathbb{Z}\colon \beta_{i,i+j}(M)\neq 0 \text{ for some } i \in \mathbb{N} \}.$
For $k \in \{0,\ldots,n\}$ we
define $d_k(M)=\min(\{j \in \mathbb{Z} \colon \beta_{k,k+j}(M)\neq 0 \} \union \{\reg(M)\})$. 
We are interested in the 
numbers $\beta_i^{k,lin}(M)=\beta_{i,i+d_k(M)}(M)$ for
$i \geq k$. Note that $\beta_i^{0,lin}(M)=\beta_i^{lin}(M)$. 
If $0 \neq \Omega_k(M)$ is the $k$-th 
syzygy module in the minimal graded free resolution of $M$
(see \cite{EI} for details), then we always have 
$\beta_{i,i+j}(M)=\beta_{i-k,i-k+j+k}(\Omega_k(M))$ for 
$i \geq k$. Therefore 
$\beta_i^{k,lin}(M)=\beta_{i-k}^{lin}(\Omega_k(M))$ for these $i$. Observe that
$d_0(\Omega_k(M))=d_k(M)+k$.

\begin{Corollary}
\label{syzkallg}
\it
Let $k \in \{0,\ldots,n \}$.
If $\beta_p^{k,lin}(M) \neq 0$ for some $p>k$, then 
$$
\beta_{p-1}^{k,lin}(M) \geq p.
$$
\end{Corollary}

For the numbers $\beta_i^{lin}(M)$ and $\beta_i^{1,lin}(M)$ we get more precise results. 
The next result was first discovered in \cite{HE}.

\begin{Theorem}
\label{syz0}
\it
Let $p \in \{0,\ldots, n\}$.
If $\beta_p^{lin}(M) \neq 0$, then 
$$
\beta_i^{lin}(M) \geq \binom{p}{i}
$$
for $i=0,\ldots, p$.
\end{Theorem}
\begin{proof}
This follows from \ref{linear} and the fact that 
there are no non-trivial boundaries in $K_i(n)_{i+d_0(M)}$.
\end{proof}

To prove lower bounds for $\beta_{i}^{1,lin}$ 
we use slightly different methods.
Let $S=K[x_1,\ldots,x_n]$.
We fix a basis $\eb$ of $L$ such that 
$\partial(e_i)=x_i$ for all $i \in [n]$ for the rest of this paper. 
For $a \in \mathbb{N}^n$ we write $x^a=x_1^{a_1} \cdots x_n^{a_n}$
and call it a monomial in $S$.
Let $F$ be a graded free $S$-module with free homogeneous basis $g_1,\ldots,g_t$.
Then
we call $x^ag_i$ a monomial in $F$ for $a \in \mathbb{N}^n$ and $i \in [t]$. 
Let $>$ be 
an arbitrary degree refining term order on $F$ with $x_1g_i>\ldots>x_ng_i>g_i$ 
(see \cite{EI} for details). 
For a homogeneous element $f \in F$ we set $\ini_{>}(f)$ for 
the maximal monomial in a presentation of $f$.
Note that we 
also defined $\ini(z)$ for some bihomogeneous $z \in K(n)$.

\begin{Lemma}
\label{trick}
\it
Let $M \subset F$ be a finitely generated graded $S$-module
and let $0 \neq z$ be a homogeneous cycle of $K_1(n,M)$ with $z=\sum_{j \geq i} m_j e_j$,
$\ini(z)=m_i e_i$. Then there exists an integer $j>i$ with $0 \neq \ini_>(m_j)>\ini_>(m_i)$.
In particular, $m_i$ and $m_j$ are $K$-linearly independent.
\end{Lemma}
\begin{proof}
We have
$$
0=\partial(z)=m_ix_i+\sum_{j>i}m_j x_j.
$$
Hence there exists an integer $j>i$
and a monomial $a_j$ of $m_j$ with 
$\ini_>(m_i)x_i=a_j x_j$ because all monomials have to cancel.
Assume that
$$\ini_>(m_i)\geq \ini_>(m_j).$$
Then
$$
\ini_>(m_i)x_i> \ini_>(m_i)x_j \geq \ini_>(m_j)x_j \geq a_j x_j
$$
is a contradiction. Therefore
$$\ini_>(m_i)< \ini_>(m_j).$$
\end{proof}

\begin{Construction}
\label{konstr}
\it
Let $M \subset F$ be a finitely generated graded $S$-module, $p \in \{0,\ldots,n\}$ 
and let $0 \neq z$ be a homogeneous cycle of $K_p(n,M)$ with $z=\sum_{J, |J|=p} m_J e_J$,
$\ini(z)=m_I e_I$. Assume that $I=\{1,\ldots,p \}$.
For $k=0,\ldots,p$ we construct inductively sets
$$J_k=\{1,\ldots,p-k,j_1,\ldots,j_k \}$$
with $j_k>p-k+1$, $j_k \neq j_i$ for $i=0,\ldots,k-1$ and
$$
\ini_>(m_{J_k})>\ini_>(m_{J_{k-1}})>\ldots>\ini_>(m_{J_0}).
$$
\rm
Set $J_0=I$.
Assume that $J_{k-1}$ is constructed. Then we apply \ref{trick} to
$$
\partial_{\{1,\ldots,p-k,j_1,\ldots,j_{k-1}\}}(z)
\text{ where } 
\ini(\partial_{\{1,\ldots,p-k,j_1,\ldots,j_{k-1}\}}(z))=m_{J_{k-1}}e_{p-k+1}
$$
and find $j_k>p-k+1$ such that 
$\ini_>(m_{\{1,\ldots,p-k,j_1,\ldots,j_{k-1}\}\union \{j_k \}})>\ini_>(m_{J_{k-1}})$.
We see that $j_k \neq j_i$ for $i=1,\ldots,k-1$ because these $e_{j_i}$ do
not appear with non-zero coefficient in $\partial_{\{1,\ldots,p-k,j_1,\ldots,j_{k-1}\}}(z)$. 
\end{Construction}

\begin{Corollary}
Let $M \subset F$ be a finitely generated graded $S$-module, 
$p \in \{0,\ldots,n \}$
and let $0 \neq z$ be a homogeneous cycle of $K_p(n,M)$ with $z=\sum_{J, |J|=p} m_J e_J$.
Then there exist $p+1$ coefficients $m_J$ of $z$, which are $K$-linearly independent.
\end{Corollary}
\begin{proof}
This follows from \ref{konstr} because the coefficients there have different leading terms.
\end{proof}

\begin{Theorem}
\label{syz1}
\it
Let $p \geq 1$ and $M$ be a finitely generated graded $S$-module with
$\beta_p^{1,lin}(M) \neq 0$. Then 
$$
\beta_i^{1,lin}(M) \geq \binom{p}{i}
$$
for $i=1,\ldots, p$.
\end{Theorem}
\begin{proof}
Let
$$
0 \to M' \to F \to M \to 0
$$
be a presentation of $M$ such that $F$ is free and $M'=\Omega_1(M)$. 
We show that
$$
\beta_i^{lin}(M') \geq \binom{p'+1}{i+1}\text{ for } p'=p-1.
$$
Since $\beta_{i}^{1,lin}(M)=\beta_{i-1}^{lin}(M')$,
this will prove the theorem.

After a suitable shift of the grading of $M$ we may assume 
that $d_0(M')=0$. Note that $M'$ is
a submodule of a free module and we can apply our construction \ref{konstr}.
Since $\beta_{p'}^{lin}(M')\neq 0$, we get a homogeneous cycle
$0 \neq z$ in $K_{p'}(n,M')_{p'}$. 
There are no boundaries except for zero in
$K_i(n,M')_i$
and therefore we only have to construct enough $K$-linearly independent
cycles in $K_i(n,M')_i$ to prove the assertion.
Assume that $\ini(z)=m_{[p']} e_{[p']}$ and construct 
the numbers $j_1,\ldots,j_{p'}$ for $z$ by \ref{konstr}.

Let $t \in \{0,\ldots,p \}$ and set $i=p'-t$. 
Consider the (by \ref{linear}) cycles $\partial_L(z)$ with 
$L \in W=W_0 \dot{\union} \ldots \dot{\union} W_t$ where 
$$
W_k=\{ I\union\{j_1,\ldots,j_k \} \colon I \subseteq [p'-k], |I|=t-k  \} \text{ for } k \in \{0,\ldots,t\}.
$$
We have
$$
|W_k|=\binom{p'-k}{t-k}
$$
and therefore
$$
|W|=\sum_{k=0}^t \binom{p'-k}{t-k}=\binom{p'+1}{t}=\binom{p'+1}{p'-t+1}=\binom{p'+1}{i+1}.
$$
If we show that the cycles $\partial_L(z)$ are $K$-linearly independent, the assertion follows.

Take $L \in W$, $L=I_L\union\{j_1,\ldots,j_{k_L} \}$ for some $I_L \subseteq [p'-k_L], |I_L|=t-k_L$.
It is easy to see by the construction \ref{konstr} that
$$
\ini(\partial_L(z))=m_{\{1,\ldots,p'-k_L \}\union\{j_1,\ldots,j_{k_L} \}}e_{\{1,\ldots,p'-k_L \}-I_L}.
$$
It is enough to show that the initial terms of the cycles are $K$-linearly independent. 

If cycles have different initial monomials in the $e_i$, there is nothing to show.
Take $L, L'$ and assume that the corresponding cycles have the same initial monomial in the $e_i$.
We have to consider two cases. If $k_L=k_{L'}$, then $I_L=I_{L'}$ and the cycles are the same.
For $k_L<k_{L'}$ the construction implies
$$
\ini(m_{\{1,\ldots,p'-k_L \}\union\{j_1,\ldots,j_{k_L}\}})<
\ini(m_{\{1,\ldots,p'-k_{L'}\}\union\{j_1,\ldots,j_{k_{L'}}\}}),
$$
which proves the $K$-linearly independence.
\end{proof}

The next corollary summarizes our results related to the conjecture 
of Herzog.
\begin{Corollary}
\label{conjecti}
\it
Let $M$ be a finitely generated graded $S$-module and $p \in \{0,\ldots,n \}$.
Suppose that $\beta_p^{lin}(M) \neq 0$ and $M$ is the $k$-th syzygy module 
in a minimal graded free resolution. Then: 
\begin{enumerate}  
\renewcommand{\labelenumi}{(\roman{enumi})}
\item
If $k=0$, then $\beta_i^{lin}(M) \geq \binom{p}{i}$ for $i=0,\ldots, p$.
\item
If $k=1$, then $\beta_i^{lin}(M) \geq \binom{p+1}{i+1}$ for $i=0,\ldots, p$.
\item
If $k>1$ and $p>0$, then $\beta_{p-1}^{lin}(M) \geq p+k$.
\end{enumerate}  
\end{Corollary}
\begin{proof}
Statement (i) was shown in \ref{syz0}.
In the proof of \ref{syz1} we proved in fact (ii). Finally, (iii) follows from 
\ref{syzkallg} since $\beta_i^{lin}(M)=\beta_{i+k}^{k,lin}(N)$ if 
$M$ is the $k$-th syzygy module in the minimal graded free resolution of some
module $N$.
\end{proof}

Recall that a finitely generated graded $S$-module $M$ satisfies 
\it Serre's condition \rm $\mathcal{S}_k$ if
$$
\depth(M_P)\geq \min(k,\dim S_P)
$$
for all $P \in \Spec(S)$.
We recall the Auslander-Bridger theorem \cite{AUBR}: 

\begin{Lemma}
\label{character}
\it
Let $M$ be a finitely generated graded $S$-module. Then $M$ is a $k$-th syzygy module 
in a graded free resolution
if and only if $M$ satisfies $\mathcal{S}_k$.
\end{Lemma}
\begin{proof}
The proof is essential the same as in \cite{EVGR} where the local case is treated.
\end{proof}

\begin{Corollary}
\label{graded}
\it
Let $M$ be a finitely generated graded $S$-module and $p \in \{0,\ldots,n \}$. 
Suppose that $M$ satisfies $\mathcal{S}_k$ and 
$\beta_p^{lin}(M) \neq 0$. 
Then:
\begin{enumerate}  
\renewcommand{\labelenumi}{(\roman{enumi})}
\item
If $k=0$, then $\beta_i^{lin}(M) \geq \binom{p}{i}$ for $i=0,\ldots, p$.
\item
If $k=1$, then $\beta_i^{lin}(M) \geq \binom{p+1}{i+1}$ for $i=0,\ldots, p$.
\item
If $k>1$ and $p>0$, then $\beta_{p-1}^{lin}(M) \geq p+k$.
\end{enumerate} 
\end{Corollary}
\begin{proof}
According to \ref{character} the module
$M$ satisfies  $\mathcal{S}_k$ if and only if $M$ is
a $k$-th syzygy module in a graded
free resolution $G_{\lpnt} \to N \to 0$ of some graded $S$-module $N$.
It is well-known (see for example \cite{AV}) that 
$G_{\lpnt}=F_{\lpnt}\dirsum H_{\lpnt}$ as graded complexes
where $F_{\lpnt}$ is the minimal graded free resolution of $N$ and
$H_{\lpnt}$ is split exact. 
Then $M$ splits as a graded module into $\Omega_k(N)\dirsum W$ 
where $W$ is a graded free $S$-module.
If $p=0$, there is nothing to show. 
For $p>0$ it follows that $\beta_p^{lin}(\Omega_k(N)) \neq 0$.
Then \ref{conjecti} applied to $\Omega_k(N)$ 
proves the corollary,
since $\beta_{i,j}(M)\geq \beta_{i,j}(\Omega_k(N))$ for all integers $i,j$.
\end{proof}

\section{Proof of the conjecture in the case of $\mathbb{Z}^n$-graded modules}
$S$ is $\mathbb{Z}^n$-graded with $\deg(x_i)=\epsilon_i=(0,\ldots,1,\ldots,0)$ 
where the $1$ is at the $i$-th position.
Let $M=\Dirsum_{a \in \mathbb{Z}^n}M_a$ be a finitely generated $\mathbb{Z}^n$-graded $S$-module.
Recall that every $\mathbb{Z}^n$-graded $S$-module $M$
is naturally $\mathbb{Z}$-graded by setting $M_i=\Dirsum_{a \in \mathbb{Z}^n, |a|=i}M_a$. 
Therefore all methods from the
last section can be applied in the following. Furthermore, the Koszul complex and homology 
are $\mathbb{Z}^n$-graded if we assign the degree $\epsilon_i$ to $e_i$. For example,
if $m \in M_a$ for some $a \in \mathbb{Z}^n$, then $\deg(m e_I)$ is $a+\sum_{i \in I} \epsilon_i$; 
or if 
$z \in K_i(j)$ is homogeneous of degree $a$, then $\deg(\partial_I(z))=a-\sum_{i \in I} \epsilon_i$.

We want to prove more precise results than in 
the last section for this more restricted situation.
Note that \ref{wichtig} and \ref{wichtig2} hold in the $\mathbb{Z}^n$-graded setting. The proofs
are verbatim the same if we replace ``graded" by ``$\mathbb{Z}^n$-graded".
We prove now a modified version of \ref{Schritt1}.
Usually we assume that $M=\Dirsum_{a \in \mathbb{N}^n}M_a$.

\begin{Lemma}
\label{Schritt1m}
\it
Let $p \in [j]$ and $t \in \mathbb{N}$.
Suppose that $H_i(j)_{i+l}=0$ for $i=p-1,\ldots,j$, $l=-1,\ldots,t-1$ and let 
$z \in K_p(j)_{p+t}$ be a $\mathbb{Z}^n$-homogeneous cycle with $0 \neq [z] \in H_p(j)_{p+j}$. 
Then there exists 
a $\mathbb{Z}^n$-homogeneous cycle $\Tilde{z}$ with:
\begin{enumerate}  
\renewcommand{\labelenumi}{(\roman{enumi})}
\item 
$[\Tilde{z}]=[z] \in H_p(j)_{p+t}$,
\item
$[\partial_{j_i}(\Tilde{z})] \in H_{p-1}(j)_{p-1+t}$ are $K$-linearly independent 
for $i=1,\ldots,p$ and some $j_i \in [j]$. 
\end{enumerate}
\end{Lemma}
\begin{proof}
We prove by induction on $j \in [n]$ that we find $[\Tilde{z}]=[z]$ and
a set $\{j_1,\ldots,j_{p} \}$ such that the cycles 
$[\partial_{j_i}(\Tilde{z})]$ are $K$-linearly independent for $i=1,\ldots,p$.

The cases $j=1$ and $j>1,p=j$ follow from \ref{wichtig} (ii) because
if $\deg(z)=a \in \mathbb{N}^n$, then all $\partial_k(z)$ have different degrees $a-\epsilon_k$
and it suffices to show that these elements are not boundaries.

Let $j>1$ and assume that $p<j$. Again it suffices to show that the cycles $\partial_{j_i}(z)$ are 
not boundaries for a suitable subset
$\{j_1,\ldots,j_{p} \}\subseteq [j]$. If such a set exists, then nothing is to prove. 
Otherwise there
exists a $k \in [j]$ with $[\partial_k(z)]=0$ and we may assume $k=j$. By \ref{wichtig}
we find $\Tilde{z}$ such that $[\Tilde{z}]=[z]$ in $H(j)$ and
$\Tilde{z}\in K(j-1)$. By \ref{homology} we can apply our induction
hypothesis and assume that $\Tilde{z}$ has the desired properties in
$H(j-1)$. Again by \ref{homology} we have 
$H_{p-1}(j-1)_{p-1+t}\subseteq H_{p-1}(j)_{p-1+t}$ and
$\Tilde{z}$ is the desired element.
\end{proof}

We need the following simple combinatorial result. Let $p \in [n]$.
Define inductively a sequence of subsets $W_i \subseteq 2^{[n]}$ for $i = 0,\ldots,p$. Set
$$
W_{0}=\{ \emptyset \}.
$$
If $W_{i-1}$ is defined, then for every set $w \in W_{i-1}$ we choose
$p-i+1$ different elements $i_1^w,\ldots, i_{p-i+1}^w $ such that $i_j^w \not\in w$.
Define
$$
W_i=\{ w \union \{i_j^w \} \colon w \in W_{i-1} \text{ and } j=1,\ldots,p-i+1 \}.
$$ 

\begin{Lemma}
\label{counting}
\it
Let $W_i$ be defined as above. Then for $i = 0,\ldots,p$ we have that
$$
|W_i|\geq \binom{p}{i}.
$$
\end{Lemma}
\begin{proof}
We prove this by induction on $p \in [n]$. The case $p=0$ is trivial, so let $p>0$ and without loss
of generality we may assume that 
$W_1=\{\{1\},\ldots,\{p\}\}$.
The set $W_i$ is the disjoint union of the sets $W_i^1=\{w \in W_i \colon 1 \in w \}$ and  
$W_i^{\hat{1}}=\{w \in W_i \colon 1 \not\in w \}$. The induction hypothesis 
applied to $W_i^1$ and $W_i^{\hat{1}}$ implies
$$
|W_i|=|W_i^1| + |W_i^{\hat{1}}|\geq \binom{p-1}{i-1} + \binom{p-1}{i}=\binom{p}{i}.
$$
\end{proof}

We prove the main theorem of this section.

\begin{Theorem}
\label{syzallg}
\it
Let $k \in [n]$ and $M$ be a finitely generated $\mathbb{Z}^n$-graded $S$-module.
If $\beta_p^{k,lin}(M) \neq 0$ for some $p \geq k$, then 
$$
\beta_i^{k,lin}(M) \geq \binom{p}{i}
$$
for $i=k,\ldots, p$.
\end{Theorem}
\begin{proof}
Without loss of generality $M=\Dirsum_{a \in \mathbb{N}^n}M_a$.
Since $\beta_p^{k,lin}(M) \neq 0$, there exists a 
$\mathbb{Z}^n$-homogeneous cycle 
$z \in K_p(n)_{p+d_k(M)}$ such that
$0 \neq[z]$ in $H_p(n)_{p+d_k(M)}$ and $\deg(z)=a$ for some $a \in \mathbb{N}^n$.
By the definition of $d_k(M)$ and \ref{homology} 
we have $H_i(n)_{i+l}=0$ for $i \geq k$ and $l=-1,\ldots,d_k(M)-1$.

We construct inductively $W_i$ as above, as well as 
cycles $z_w$ for each $w \in W_i$ such that $[z_w]\neq 0$, 
$\deg(z_w)=a-\sum_{i \in w} \epsilon_i$,  $0 \neq [\partial_{i^w_k}(z_w)]$ for $k=1,\ldots,p-i$ and 
suitable $i^w_k \not\in w$. 
Furthermore, $z_w$ is an element of the Koszul complex with respect to the variables $x_i$ 
with $i \not\in w$. 
For $i \geq k$ we take all cycles $z_w \in K_i(n)_{i+d_k(M)}$ 
with $w \in W_i$ which have different $\mathbb{Z}^n$-degree. 
They are not zero and therefore
$K$-linearly independent in homology. 
By \ref{counting} there are at least $\binom{p}{p-i}=\binom{p}{i}$ 
of them and this concludes the proof.

Let $W_0=\{ \emptyset \}$.
By \ref{Schritt1m} we can choose $z$ in a way such that $[z_{i_j}]=[\partial_{i_j}(z)]\neq 0$ for
$j=1,\ldots,p$ and some $i_j \in [n]$.
Choose $z_{\emptyset}=z$ and $i^{\emptyset}_j=i_j$.

If $W_{i-1}$ and $z_w$ for $w \in W_{i-1}$ are constructed, then 
define $W_i$ with $W_{i-1}$ and the given $i^w_k$ for $w \in W_{i-1}$.
For $w' \in W_i$ with $w'=w \union \{ i^w_k \}$, re-choose   
$z_{w'}=\partial_{i^w_k}(z_w)$ by \ref{Schritt1m} in such a way that 
$[\partial_{i_j^{w'}}(z_{w'})]\neq 0$ for $j=1,\ldots,p-i$ and some $i_j^{w'} \in [n]$.

Note that since $z_{w}$ has no monomial which is divided by some $e_i$ for $i \in w$, we
can use \ref{wichtig} and \ref{wichtig2} to avoid these $e_i$ in the
construction of $z_{w'}$ again. 
By \ref{homology} the cycles $\partial_{i_j^{w'}}(z_{w'})$ are also not zero in $H(n)$.
Clearly $i_j^{w'} \not\in w'$ and the assertion follows.
\end{proof}

In the $\mathbb{Z}^n$-graded setting we prove the desired results about 
$\beta_i^{lin}$ in full
generality.
\begin{Corollary}
\label{kth}
\it
Let $M$ be a finitely generated $\mathbb{Z}^n$-graded $S$-module 
and $M$ is the $k$-th syzygy module 
in a minimal $\mathbb{Z}^n$-graded free resolution. 
If $\beta_p^{lin}(M) \neq 0$ for some $p \in \mathbb{N}$, then 
$$
\beta_i^{lin}(M) \geq \binom{p+k}{i+k}
$$
for $i=0,\ldots, p$.
\end{Corollary}
\begin{proof}
This follows from \ref{syzallg} and the fact that
$$
\beta_i^{lin}(M)=\beta_{i+k}^{k,lin}(N)\geq \binom{p+k}{i+k}
$$
where $M$ is the $k$-th syzygy module of a $\mathbb{Z}^n$-graded $S$-module $N$. 
\end{proof}

Analogue to \ref{character} we get:
\begin{Lemma}
\label{Sk}
\it
Let $M$ be a finitely generated $\mathbb{Z}^n$-graded $S$-module.
Then $M$ satisfies $\mathcal{S}_k$ if and only if $M$ is a $k$-th syzygy module
in a $\mathbb{Z}^n$-graded free resolution.
\end{Lemma}

\begin{Corollary}
\label{}
\it
Let $M$ be a finitely generated $\mathbb{Z}^n$-graded $S$-module and $M$ satisfies $\mathcal{S}_k$.
If $\beta_p^{lin}(M) \neq 0$ for some $p \in \mathbb{N}$, then 
$$
\beta_i^{lin}(M) \geq \binom{p+k}{i+k}
$$
for $i=0,\ldots, p$.
\end{Corollary}
\begin{proof}
The assertion follows from \ref{kth}
with similar arguments as in the graded case.
\end{proof}

\section{Bounds for Betti numbers of ideals with a fixed number of generators
in given degree and a linear resolution}
In this section we are interested in bounds for the graded Betti numbers of
graded ideals of $S$. We assume that
the field $K$ is infinite and fix a basis $\xb=x_1,\ldots,x_n$ of $S_1$. 
For a monomial $x^a$ of $S$ with $|a|=d$ we denote
by $R(x^a)=\{x^b \colon |b|=d, x^b \geq_{rlex} x^a \}$ the revlex-segment of $x^a$ 
where
$>_{rlex}$ is defined as follows:
$x^a >_{rlex}x^b$ if either
$|a|>|b|$ or
$|a|=|b|$ and there exists an integer $r$ such that
$a_r<b_r$ and $a_s=b_s$ for $s>r$.
Note that for a given $d \in \mathbb{N}$ and $0 \leq k \leq \binom{n+d-1}{d}$ there exists 
a unique ideal $I(d,k)$ which is generated in degree $d$ by a revlex-segment
$R(x^a)$ for some monomial $x^a$ with $|R(x^a)|=k$.

Following Eliahou and Kervaire \cite{ELKE} 
for a given monomial $x^a$ let $m(x^a)$ be the maximal $i$ with $x_i$ divides $x^a$. 
An ideal is a monomial ideal if it is generated by monomials.
We call a monomial ideal $I$ stable if
for all monomials $x^a \in I$ we have $x_ix^a/x_{m(x^a)} \in I$ for $i<m(x^a)$. 
Is it easy to see that it is enough to prove this condition
for the generators of the ideal $I$. For example
$I(d,k)$ is stable.
   
For stable ideals there exist explicit formulas for the Betti numbers (see \cite{ELKE}). 
Let $I \subset S$ be a stable ideal and $G(I)$ be the set of minimal generators of $I$. Then
$$
\beta_{i,i+j}(I)=\sum_{x^a \in G(I), |a|=j}\binom{m(x^a)-1}{i} \quad(*).
$$
If a stable ideal $I$ is generated in one degree $d$, then $I$ has a linear resolution.

\begin{Proposition}
\label{stable}
Let $I \subset S$ be a stable ideal generated in degree $d \in \mathbb{N}$ with 
$\beta_{0,d}(I)=k$. Then 
$$
\beta_{i,i+j}(I) \geq \beta_{i,i+j}(I(d,k))
$$
for all $i,j \in \mathbb{N}$.
\end{Proposition}
\begin{proof}
Fix $d \in \mathbb{N}$ and $0 \leq k \leq \binom{n+d-1}{d}$.
Let $l(I)$ be the number of monomials in $I_d$ which are not monomials in $I(d,k)_d$.
We prove the statement by induction on $l(I)$.
If $l(I)=0$, then $I=I(d,k)$ and there is nothing to show.

Assume that $l(I)>0$. Let $x^a$ be the smallest monomial in $I_d$ with respect to
$>_{rlex}$ which is not in $I(d,k)_d$, and $x^b$ be the largest monomial which is
in $I(d,k)_d$, but not in $I_d$.
Define the ideal $\Tilde{I}$ by $G(\Tilde{I})=(G(I)\setminus \{x^a \})\union\{x^b \}$.
Then $l(\Tilde{I})=l(I)-1$, and $\Tilde{I}$ is also stable.
Thus, by the induction hypothesis
$$
\beta_{i,i+j}(\Tilde{I}) \geq \beta_{i,i+j}(I(d,k)).
$$
Since $x^b>_{rlex} x^a$, the revlex order implies
$
m(x^b)\leq m(x^a).$
Therefore $(*)$ yields
$$
\beta_{i,i+j}(I) \geq \beta_{i,i+j}(\Tilde{I}),
$$
which proves the assertion.
\end{proof}   

For a graded ideal $I$, the initial ideal $\ini(I)$ is generated
by all $\ini(f)$ for $f \in I$ with respect to some monomial order. 

Every element $g$ of the general linear group $\GL(n)$ induces a linear automorphism of $S$ by
$$
g(x_j)=\sum_{i=1}^{n} g_{i,j}x_i \text{ for } g=(g_{i,j}).
$$
There is a non-empty open
set $U \subset \GL(n)$ and a unique monomial ideal $J$ with 
$J=\ini(g(I))$ for every $g \in U$ with respect to the revlex order (for details see \cite{EI}).
We call $J$ the generic initial ideal of $I$ and denote it by $\Gin(I)$.
A nice property is that $\Gin(I)$ is Borel-fixed, i.e.
$\Gin(I)=b \Gin(I)$ for all $b \in B$ where $B$ is the Borel subgroup of $\GL(n)$ which
is generated by all upper triangular matrices.

\begin{Proposition}
\label{borel}
Let $d \in \mathbb{N}$ and $I \subset S$ be a graded ideal with $d$-linear resolution. Then
$\Gin(I)$ is stable, independent of the characteristic of $K$, and
$\beta_{i,i+j}(I)=\beta_{i,i+j}(\Gin(I))$ for all $i,j \in \mathbb{N}$.
\end{Proposition}
\begin{proof}
It is well-known that $\reg(\Gin(I))=\reg(I)$. Therefore
$\reg(\Gin(I))=d$ and $\Gin(I)$ also has  a $d$-linear resolution.
\cite[Prop. 10]{EIRETO} implies that
a Borel-fixed monomial ideal, which is generated in degree $d$, 
has regularity $d$ if and only if it is stable. Thus
we get that $\Gin(I)$ is a stable ideal, independent of the characteristic of $K$.

Since $I$ has a linear resolution, we obtain by the main result in \cite{ARHEHI} that
$\beta_{i,i+j}(I)=\beta_{i,i+j}(\Gin(I))$ for all $i,j \in \mathbb{N}$.
\end{proof}

\begin{Theorem}
Let $d \in \mathbb{N}$, $0 \leq k \leq \binom{n+d-1}{d}$ 
and $I \subset S$ be a graded ideal with $d$-linear resolution and $k$ generators.
Then 
$$
\beta_{i,i+j}(I)\geq \beta_{i,i+j}(I(d,k))
$$
for all $i,j \in \mathbb{N}$.
\end{Theorem}
\begin{proof}
This follows from \ref{stable} and \ref{borel}.
\end{proof}

Consider the lexicographic order $>_{lex}$. Recall
that for monomials $x^a,x^b \in S$ we have 
$x^a >_{lex}x^b$ if either
$|a|>|b|$ or
$|a|=|b|$ and there exists an integer $r$ such that
$a_r>b_r$ and $a_s=b_s$ for $s<r$.
Then
$L(x^a)=\{x^b \colon |b|=d, x^b \geq_{lex} x^a \}$ is the lex-segment of $x^a$.
For a given $d \in \mathbb{N}$ and $0 \leq k \leq \binom{n+d-1}{d}$ there exists 
a unique ideal $J(d,k)$ which is generated in degree $d$ by a lex-segment
$L(x^a)$ for some monomial $x^a$ with $|L(x^a)|=k$.
It is easy to see that $J(d,k)$ is a lex-ideal, i.e.\  
if $x^b>_{lex}x^c$ and $x^c \in J(d,k)$, then $x^b \in J(d,k)$. 
In particular $J(d,k)$ is stable.

\begin{Proposition}
\label{stable2}
Let $d \in \mathbb{N}$, $0 \leq k \leq \binom{n+d-1}{d}$ 
and $I \subset S$ be a graded ideal with $d$-linear resolution and $k$ generators.
Then 
$$
\beta_{i,i+j}(I)\leq \beta_{i,i+j}(J(d,k))
$$
for all $i,j \in \mathbb{N}$.
\end{Proposition}
\begin{proof}
By \cite[Thm. 31]{PA} we find a lex-ideal $L$ with the same
Hilbert function as $I$ and  
$
\beta_{i,i+j}(I)\leq\beta_{i,i+j}(L)
$
for all $i,j \in \mathbb{N}$.
We see that
$
\beta_{0,d}(L)=\beta_{0,d}(I)=k
$
because these ideals share the same Hilbert function.
It follows that
$J(d,k)=(L_d)$, and in particular $G(J(d,k))=G(L)_d$.
Therefore
$$
\beta_{i,i+d}(I)\leq\beta_{i,i+d}(L)=\beta_{i,i+d}(J(d,k))
$$
for all $i,j \in \mathbb{N}$ where the last equality follows from $(*)$.
\end{proof}

Fix $d \in \mathbb{N}$ and
$0 \leq k \leq\binom{n+d-1}{d}$.
Let $\mathcal{B}(d,k)$ be the set of 
Betti sequences $\{\beta_{i,j}(I)\}$ where $I$ is a
graded ideal with $d$-linear resolution and 
$\beta_{0,d}(I)=k$. 
On $\mathcal{B}(d,k)$ we consider a partial 
order: We set $\{\beta_{i,j}(I)\}\geq \{\beta_{i,j}(J)\}$ 
if $\beta_{i,j}(I)\geq \beta_{i,j}(J)$ 
for all $i,j \in \mathbb{N}$.
\begin{Corollary}
Let $d \in \mathbb{N}$ and
$0 \leq k \leq\binom{n+d-1}{d}$.
Then $\{\beta_{i,j}(I(d,k))\}$ is the unique minimal element and
$\{\beta_{i,j}(J(d,k))\}$ is the unique maximal element of $\mathcal{B}(d,k)$.
\end{Corollary}

\newpage
\noindent
Tim R\"omer\\
FB6 Mathematik und Informatik\\
Universit\"at Essen\\
45117 Essen\\
Germany\\
tim.roemer@uni-essen.de


\begin{thebibliography}{99}
\bibitem{ARHEHI}
A. Aramova, J. Herzog and T. Hibi, 
Ideals with stable Betti numbers, 
\it Adv. Math. \rm 
\bf 152 \rm (2000),
no. 1,
72-77.
\bibitem{AUBR}
M. Auslander and M. Bridger,
Stable Module Theory, 
\it Mem. Amer. Math. Soc. \rm (1969),no. 94.
\bibitem{AV}
L. L. Avramov,
Infinite free resolutions,
Six lectures on commutative algebra (Bellaterra, 1996), 
1-118,
\it Progr. Math. \rm \bf 166\rm, 
Birkh\"auser, 
Basel, 
1998. 
\bibitem{BI}
A. Bigatti, 
Upper bounds for the Betti numbers of a given Hilbert function, 
\it Comm. Algebra \rm 
\bf 21 \rm (1993),
no. 7,  
2317-2334.
\bibitem{BRHE}
W. Bruns and J. Herzog, 
Cohen-Macaulay rings, 
revised edition,
\it Cambridge Studies in Advanced Mathematics \rm \bf 39\rm, 
Cambridge Univ. Press,
Cambridge,
1998.
\bibitem{CHEV}
H. Charalambous and E. Evans, 
Resolutions with a given Hilbert function, 
\it Contemp. Math. \rm 
\bf 159 \rm (1994),
19-26.
\bibitem{EI}
D. Eisenbud, 
Commutative algebra with a view toward algebraic geometry,
\it Graduate Texts in Mathematics \rm \bf 150\rm, 
Springer-Verlag,
New York,
1995.
\bibitem{EIKO}
D. Eisenbud and J. Koh, 
Some linear syzygy conjectures, 
\it Adv. Math. \rm 
\bf 90 \rm (1991), 
no. 1,
47-76.
\bibitem{EIRETO}
D. Eisenbud, A. Reeves and B. Totaro, 
Initial ideals, Veronese subrings, and rates of algebras, 
\it Adv. Math. \rm 
\bf 109 \rm (1994), 
no. 2,
168-187. 
\bibitem{ELKE}
S. Eliahou and M. Kervaire, 
Minimal resolutions of some monomial ideals,
\it J. Algebra \rm
\bf 129 \rm (1990), 
no. 1,
1-25.
\bibitem{EVGR}
E. G. Evans and P. Griffith, 
Syzygies, 
London Mathematical Society, 
\it Lecture Note Series \rm \bf 106\rm, 
Cambridge Univ. Press, 
Cambridge-New York, 
1985.  
\bibitem{GR}
M. Green, 
Koszul cohomology and the geometry of projective varieties, 
\it J. Diff. Geom. \rm 
\bf 19 \rm (1984), 
no. 1,
125-171.
\bibitem{HE}
J. Herzog, 
The linear strand of a graded free resolution, 
unpublished notes (1998).
\bibitem{HU}
H. Hulett, 
Maximum Betti numbers of homogeneous ideals 
with a given Hilbert function, 
\it Comm. Algebra \rm 
\bf 21 \rm (1993), 
no. 7,
2335-2350.
\bibitem{PA}
K. Pardue, 
Deformation classes of graded modules and maximal Betti numbers, 
\it Illinois J. Math. \rm 
\bf 40 \rm (1996), 
no. 4,
564-585.
\bibitem{REWE}
V. Reiner and V. Welker, 
Linear syzygies of Stanley-Reisner ideals,
\it Math. Scand., \rm 
to appear.
\end{thebibliography}
\end{document}